\documentclass[10pt]{article}

\usepackage{amssymb,latexsym}
\usepackage{amsfonts}
\usepackage{amsmath,amsthm,graphicx}

\newcounter{cnt1}
\newcounter{cnt2}
\newcounter{cnt3}
\newcommand{\blr}{\begin{list}{$($\roman{cnt1}$)$} {\usecounter{cnt1}
        \setlength{\topsep}{0pt} \setlength{\itemsep}{0pt}}}
\newcommand{\bla}{\begin{list}{$($\alph{cnt2}$)$} {\usecounter{cnt2}
        \setlength{\topsep}{0pt} \setlength{\itemsep}{0pt}}}
\newcommand{\bln}{\begin{list}{$($\arabic{cnt3}$)$} {\usecounter{cnt3}
                \setlength{\topsep}{0pt} \setlength{\itemsep}{0pt}}}
\newcommand{\el}{\end{list}}

\newtheorem{thm}{Theorem}[section]

\newtheorem{prop}[thm]{Proposition}
\newtheorem{defn}[thm]{Definition}

\newtheorem{rem}[thm]{Remark}

\newcommand{\ilim}{\mathop{\varprojlim}\limits}

\begin{document}

\title{Leibniz algebra deformations of a Lie algebra}
\author{Alice Fialowski
\thanks{The work was partially supported by grants from INSA India and HAS Hungary.}~  and Ashis Mandal
}
\maketitle
\date{}
\begin{abstract}
In this note we compute Leibniz algebra deformations of the $3$-dimensional nilpotent Lie algebra $\mathfrak{n}_3$ and compare it with its Lie deformations. It turns out that there are $3$ extra Leibniz deformations. We also describe the versal Leibniz deformation of $\mathfrak{n}_3$ with the versal base. 
\end{abstract}

%\footnote{AMS Mathematics Subject Classification : 54H20, 57S25.}
%\newsymbol \rtimes 226F

\section{Introduction}
Since a Lie algebra is also a Leibniz algebra, a natural question arises. If we consider a Lie algebra as a Leibniz algebra and compute its Leibniz algebra deformations, is it true that we can get more Leibniz algebra deformations, than just the Lie deformations of the original Lie algebra ?

In this note we will demonstrate the problem on a three dimensional Lie algebra example for which we completely describe its versal Lie deformation and versal Leibniz deformation. It turns out that beside the Lie deformations we get three non-equivalent Leibniz deformations which are not Lie algebras. 

Our example is the following.
Consider a three dimensional vector space $L$ spanned by  $\{e_1,~e_2,~e_3\}$ over $\mathbb{C}$. Define the Heisenberg Lie algebra $\mathfrak{n}_3$ on it with the bracket matrix
$$A= \left( \begin{array}{lll}
0& 0 & 1\\
0& 0 & 0 \\
0& 0 & 0  \end{array}\right).$$ 
Here the columns are the Lie brackets $[e_1,e_2],[e_1,e_3]$ and $[e_2,e_3]$.
That means $[e_2, e_3]=e_1$  and all the other brackets are zero (except of course  $[e_3,e_2]=-e_1$).
This is the only nilpotent three dimensional Lie algebra. We compute infinitesimal Lie and Leibniz deformations and show that there are three additional Leibniz cocycles beside the five Lie cocycles. We will show that all these infinitesimal deformations are extendable without any obstructions. We also describe the versal Leibniz deformation.

The structure of the paper is as follows. In Section $1$ we recall the necessary preliminaries about Lie and Leibniz cohomology and deformations. In Section $2$ we recall the classification of three dimensional Lie algebras and describe all non-equivalent deformations of the nilpotent Lie algebra $\mathfrak{n}_3$ . In Section $3$ we give the classification of three dimensional nilpotent Leibniz algebras. Among those $\mathfrak{n}_3$ is the only nontrivial Lie algebra. Then we compute Leibniz cohomology and give explicitly all non-equivalent infinitesimal deformations. In Section $4$ we show that all infinitesimal deformations are extendable as all the Massey squares turn out to be zero. We identify our deformations with the classified objects. Finally we  show that the versal Leibniz deformation is the universal infinitesimal one, and describe the base of the versal deformation.

\section{Preliminaries }
Let us recall first the Lie algebra cohomology.
\begin{defn}
Suppose $\mathfrak{g}$ is a Lie algebra and $A$ is a module over $\mathfrak{g}$. Then a $q$-dimensional cochain of the Lie algebra $\mathfrak{g}$ with coefficients in $A$ is a skew -symmetric $q$-linear map on $\mathfrak{g}$ with values in $A$; the space of all such cochains is denoted by $C^q(\mathfrak{g};A)$. Thus, $C^q(\mathfrak{g};A)=Hom(\Lambda^q \mathfrak{g},A)$; this last representation transforms $C^q(\mathfrak{g};A)$ into a $\mathfrak{g}$-module. The differential $$d=d_q:C^q(\mathfrak{g};A)\longrightarrow C^{q+1}(\mathfrak{g};A)$$
is defined by the formula
\begin{equation*}
 \begin{split}
dc(g_1,\cdots,g_{q+1})=&\sum_{1\leq s<t\leq q+1}(-1)^{s+t-1}c([g_s,g_t],g_1,\cdots,\hat{g}_s,\cdots,\hat{g}_t,\cdots,g_{q+1})\\
&\sum_{1\leq s\leq q+1 }(-1)^{s}[g_s,c(g_1,\cdots,\hat{g}_s,\cdots,g_{q+1})],
 \end{split}
\end{equation*}
where $c \in C^q(\mathfrak{g};A)$ and $g_1,\cdots,g_{q+1} \in \mathfrak{g}$.

We complete the definition by putting $C^q(\mathfrak{g};A)=0$ for $q<0$, $d_q=0$ for $q<0$. As can be easily checked, $d_{q+1}\circ d_q=0$ for all $q$, so that $(C^q(\mathfrak{g},d_q))$ is an algebraic complex; this complex is denoted by $C^{*}((\mathfrak{g};A))$, while the corresponding cohomology is referred to as the cohomology of the Lie algebra $\mathfrak{g}$ with coefficients in $A$ and is denoted by $H^q(\mathfrak{g};A)$.
\end{defn}
Leibniz algebras were introduced by J.-L. Loday \cite{L1,LP}. Let $\mathbb{K}$ denote a field.
\begin{defn}
A Leibniz algebra is a $\mathbb{K}$-module $L$, equipped with a bracket operation that satisfies the Leibniz identity: 
$$[x,[y,z]]= [[x,y],z]-[[x,z],y],~~\mbox{for}~x,~y,~z \in L.$$ 
\end{defn}

Any Lie algebra is automatically a Leibniz algebra, as in the presence of antisymmetry, the Jacobi identity is equivalent  to the Leibniz identity. More  examples of Leibniz algebras were given in \cite {L1,L2,LP}, and recently for instance in \cite{A3,A2}.

Let $L$ be a Leibniz algebra and $M$ a representation of $L$. By definition, $M$ is a $\mathbb{K}$-module equipped with two actions (left and right) of $L$,
$$[-,-]:L\times M\longrightarrow M~~\mbox{and}~[-,-]:M \times L \longrightarrow M ~~\mbox{such that}~$$ 
$$[x,[y,z]]=[[x,y],z]-[[x,z],y]$$
holds, whenever one of the variables is from $M$ and the two others from $L$.

Define $CL^n({L}; {M}):= \mbox{Hom} _\mathbb{K}({L}^{\otimes n}, {M}), ~n\geq 0.$ Let 
$$\delta^n : CL^n({L}; {M})\longrightarrow CL^{n+1}(L; M)$$ 
be a $\mathbb{K}$-homomorphism defined by 
\begin{equation*}
\begin{split}
&\delta^nf(x_1,\cdots,x_{n+1})\\
&:= [x_1,f(x_2,\cdots,x_{n+1})] + \sum_{i=2}^{n+1}(-1)^i[f(x_1,\cdots,\hat{x}_i,\cdots,x_{n+1}),x_i]\\
&+ \sum_{1\leq i<j\leq n+1}(-1)^{j+1}f(x_1,\cdots,x_{i-1},[x_i,x_j],x_{i+1},\cdots,\hat{x}_j,\cdots, x_{n+1}).
\end{split}
\end{equation*}
Then $(CL^*(L; M),\delta)$ is a cochain complex, whose cohomology is called the cohomology of the Leibniz algebra $L$ with coefficients in the representation $M$. The $n$th cohomology is denoted by $HL^n(L; M)$. In particular, $L$ is a representation of itself with the obvious action given by the bracket in $L$. The  $n$th cohomology of $L$ with coefficients in itself is denoted by $HL^n(L; L).$

Let $S_n$ be the symmetric group of $n$ symbols. Recall that a permutation $\sigma \in S_{p+q}$ is called a $(p,q)$-shuffle, if $\sigma(1)<\sigma(2)<\cdots<\sigma(p)$, and $\sigma(p+1)<\sigma(p+2)<\cdots<\sigma(p+q)$. We denote the set of all $(p,q)$-shuffles in $S_{p+q}$ by $Sh(p,q)$.

For $\alpha \in CL^{p+1}(L;L)$ and $\beta \in CL^{q+1}(L;L)$, define $\alpha \circ \beta \in CL^{p+q+1}(L;L)$ by
\begin{equation*}
\begin{split}
&\alpha \circ \beta (x_1,\ldots,x_{p+q+1} )\\
=&~\sum_{k=1}^{p+1}(-1)^{q(k-1)}\{\sum_{\sigma \in Sh(q,p-k+1)}sgn(\sigma)\alpha(x_1,\ldots,x_{k-1},\beta(x_k,x_{\sigma(k+1)},\ldots,x_{\sigma(k+q)}),\\
&~~~~~~~~~~~~~~~~~~~~~~~~~~~~~~~~~~~~~~~~~~~~~~~~~~~~~x_{\sigma(k+q+1)},\ldots,x_{\sigma(p+q+1)}) \}.
\end{split}
\end{equation*}

The graded cochain module 
$CL^{*}(L;L)=\bigoplus_{p} CL^p(L;L)$ equipped with the bracket defined by  
$$[\alpha,\beta]=\alpha \circ \beta + (-1)^{pq+1} \beta \circ \alpha
~~\mbox{for}~ \alpha \in CL^{p+1}(L;L)~~\mbox{and}~\beta \in CL^{q+1}(L;L)$$
 and the differential map $d$  by $d \alpha =(-1)^{|\alpha|}\delta \alpha~\mbox{for}~\alpha \in CL^{*}(L;L) $ is a differential graded Lie algebra.

Let now $\mathbb{K}$ a field of zero characteristic and the tensor product over $\mathbb{K}$ will be denoted by $\otimes$. We recall the notion of deformation of a Lie (Leibniz) algebra $\mathfrak{g}$ ($L$) over a commutative algebra base $A$ with a fixed augmentation $\varepsilon:{A}\rightarrow
\mathbb{K}$ and maximal ideal $\mathfrak{M}$. Assume $dim(\mathfrak{M}^k/\mathfrak{M}^{k+1})<\infty$ for every $k$ (see \cite{FMM}).
\begin{defn}
A deformation $\lambda$ of a Lie algebra $\mathfrak{g}$ (or a Leibniz algebra ${L}$) with base
$({A},\mathfrak{M})$, or simply with base ${A}$ is an $A$-Lie algebra (or an ${A}$-Leibniz
algebra) structure on the tensor product
${A}\otimes {\mathfrak {g}}$ (or $A\otimes L$) with the bracket $[,]_\lambda$ such that
 \[
 \varepsilon\otimes id:{A}\otimes {\mathfrak{g}}\rightarrow \mathbb{K}\otimes {\mathfrak{g}}~~ (\mbox{or}~\varepsilon\otimes id:{A}\otimes {L}\rightarrow \mathbb{K}\otimes {L})
 \]
 is an $A$-Lie algebra (${A}$-Leibniz algebra)  homomorphism. 
\end{defn}
A deformation of the Lie (Leibniz) algebra $\mathfrak{g}$ ($L$) with base $A$ is called {\it infinitesimal}, {or \it first order}, if in addition to this $\mathfrak{M}^2=0$. We call a deformation of {\it order k}, if $\mathfrak{M}^{k+1}=0$. A deformation with base is called local if $A$ is a local algebra over $\mathbb{K}$, which means $A$ has a unique maximal ideal.

Suppose $A$ is a complete local algebra ( $A=\ilim_{n\rightarrow
\infty}({A}/{\mathfrak{M}^n})$), where $\mathfrak{M}$ is the maximal
ideal in $A$. Then a deformation of $\mathfrak{g}$ ($L$) with base $A$ which is obtained as the projective limit of deformations of $\mathfrak{g}$ $(L)$ with base $A/\mathfrak{M}^{n}$ is called a {\it formal deformation} of $\mathfrak{g}$ $(L)$.

\begin{defn}
Suppose $\lambda$ is a given deformation of $L$ with base
$(A,\mathfrak{M})$ and augmentation $\varepsilon:{A}\rightarrow
\mathbb{K}$. Let $A^\prime$ be another commutative algebra with
identity and a fixed augmentation
$\varepsilon^{\prime}:{A^\prime}\rightarrow \mathbb{K}$. Suppose
$\phi:A \rightarrow A^{\prime} $ is an algebra homomorphism with
$\phi(1)=1$ and $\varepsilon^{\prime} \circ \phi =\varepsilon$. Let
$ker(\varepsilon^{\prime})= \mathfrak{M}^\prime$. Then the push-out
$\bf{\phi_{*} \lambda}$ is the deformation of $L$ with base
$(A^\prime,\mathfrak{M}^\prime)$ and bracket
          $$[{a_1 }^\prime \otimes_A (a_1\otimes {l_1}),a_2 ^\prime
\otimes_A(a_2\otimes l_2) ]_{\phi_* \lambda}=a_1 ^\prime a_2 ^\prime
\otimes_A[a_1\otimes l_1,a_2\otimes l_2]_\lambda $$
 where $a_1^\prime,a_2 ^\prime \in {A}^\prime,~ a_1,a_2 \in A$ and
$l_1,l_2 \in L$. Here $A^\prime$ is considered as an $A$-module by the
map $a^\prime \cdot a=a^\prime \phi(a)$ so that $$A^\prime \otimes
L=(A^\prime {\otimes}_{A} A)\otimes L =A^\prime {\otimes}_{A}(A \otimes
L).$$
\end{defn}

\begin{defn}(see \cite{F}) 
Let $C$ be a complete local algebra. A formal deformation $\eta$ of a Lie algebra $\mathfrak{g}$ (Leibniz algebra $L$) with base $C$ is called versal, if\\
(i)~for any formal deformation $\lambda$ of $\mathfrak{g}$ ($L$) with  base $A$ there exists a homomorphism $f:C \rightarrow A$ such that the deformation $\lambda$ is equivalent to $f_{*}\eta$; \\
(ii)~if $A$ satisfies the condition ${\mathfrak{M}}^2=0$, then $f$ is unique. 
\end{defn} 

\begin{thm}
If $H^2(\mathfrak{g};\mathfrak{g})$ is finite dimensional, then there exists a versal deformation of $\mathfrak{g}$ (similarly for $L$).
\end{thm}
\begin{proof}
Follows from the general theorem of Schlessinger \cite{Sch}, like it was shown for Lie algebras in \cite{F}.
\end{proof}

In \cite{FiFu} a construction for a versal deformation of a Lie algebra was given and it was generalized to Leibniz algebras in \cite{FMM}. The computation of a  specific example is given in \cite{AM2}.

Let us describe the universal infinitesimal deformation (see \cite{FiFu} and \cite{FMM}).
For simplicity we will only discuss the Leibniz algebra case.
 
Assume that $dim (HL^2(L;L)) < \infty$. Denote the space $HL^2(L;L)$ by $\mathbb{H}$. Consider the algebra $C_1=\mathbb{K}\oplus \mathbb{H}^\prime$  where $\mathbb{H}^\prime$ is the dual of $\mathbb{H}$~, by setting $$(k_1,h_1)\cdot(k_2,h_2)=(k_1 k_2,k_1 h_2+k_2 h_1)~\mbox{for}~(k_1,h_1), (k_2,h_2)\in C_1.$$   Observe that the second summand is an ideal of $C_1$ with zero multiplication. Fix a homomorphism $$\mu: \mathbb{H} \longrightarrow
CL^2(L;L)=Hom(L^{\otimes 2};L) $$
which takes a cohomology class into a cocycle representing it. Notice that there is an isomorphism $\mathbb{H}^\prime \otimes L \cong Hom(\mathbb{H}~;L)$, 
so we have
$$C_1 \otimes L 
 =(\mathbb{K} \oplus \mathbb{H}^\prime)\otimes L 
 \cong (\mathbb{K}\otimes L) \oplus (\mathbb{H}^\prime \otimes L) 
 \cong L \oplus Hom(\mathbb{H}~;L).$$
Using the above identification, define a Leibniz bracket on $C_1 \otimes L$ as follows.
For $(l_1,\phi_1),(l_2,\phi_2) \in L \oplus Hom(\mathbb{H}~;L)$ let
$$[(l_1,\phi_1),(l_2,\phi_2)]_{\eta_1}=([l_1,l_2],\psi)$$ where the map  $\psi:\mathbb{H} \longrightarrow L$ is given by 
$$\psi(\alpha)=\mu(\alpha)(l_1,l_2)+[\phi_{1}(\alpha),l_2]+[l_1,\phi_2(\alpha)]~\mbox{for}~\alpha \in \mathbb{H}~.$$ 
It is straightforward to check that $C_1\otimes L$ along with the above bracket $\eta_1$ is a Leibniz algebra over $C_1$. The Leibniz identity is a consequence of the fact that $\delta \mu(\alpha)=0~ \mbox{for}~ \alpha \in \mathbb{H}$~. Thus $\eta_1$ is an infinitesimal deformation of $L$ with base $C_1=\mathbb{K}\oplus \mathbb{H}^\prime $. It is proved in \cite{FMM}:
\begin{prop}\label{up toisomorphism}
Up to an isomorphism, the deformation $\eta_1$ does not depend on
the choice of $\mu$. 
\end{prop} 
  
\begin{rem}\label{exp of inf}
Suppose $\{h_i \}_{1\leq i \leq n}$ is a basis of $\mathbb{H}$ and $\{g_i\}_{1\leq i \leq n}$ is the dual basis. Let $\mu(h_i)=\mu_i \in CL^2(L;L)$. Under the identification $C_1 \otimes L = L \oplus Hom(\mathbb{H}~;L)$, an element $(l,\phi)\in L \oplus Hom(\mathbb{H}~;L)$ corresponds to $1\otimes l +\sum_{i=1}^{n}{g_i\otimes \phi(h_i)}$. 
Then for $(l_1,\phi_1),(l_2,\phi_2) \in L \oplus Hom(\mathbb{H};L)$ their bracket $([l_1,l_2],\psi)$ 
corresponds to 
$$1\otimes [l_1,l_2]+ \sum_{i=1}^{n} g_i\otimes (\mu_i(l_1,l_2)+[\phi_1(h_i),l_2]+[l_1,\phi_2(h_i)]).$$
In particular, for $l_1,l_2 \in L$ we have 
$$[1\otimes l_1,1\otimes l_2]_{\eta_1}=1\otimes [l_1,l_2]+\sum_{i=1}^{n}g_i \otimes \mu_i(l_1,l_2).$$ 
\end{rem}
The main property of $\eta_{1}$ is the universality in the class of infinitesimal deformations with a finite dimensional base. 
 
\begin{prop}\label{couniversal}
For any infinitesimal deformation $\lambda$ of a Leibniz algebra
$L$ with a finite dimensional  base $A$ there exists a unique
homomorphism $\phi:C_1=({\mathbb{K}}\oplus
\mathbb{H}^\prime)\longrightarrow A$ such that $\lambda$ is
equivalent to the push-out $\phi_{*}\eta_1$.
\end{prop}
After obtaining the universal infinitesimal deformation we could like to extend it to higher order. If we have a universal infinitesimal deformation with basis cocycles  $\{\phi_i\}_{i=i}^r$, the obstructions to extend it to a second order deformation are the Lie brackets $[\phi_i,\phi_j]$ in the cochain complex (see \cite{F,FMM}). This are also called first order Massey operations. If these bracket cochains are coboundaries we can extend our infinitesimal to the second order deformation. The construction of a versal deformation for Leibniz algebras is given in \cite{FMM}.

\section{Deformations of the three dimensional Lie algebra $\mathfrak{n}_3$}
Let us recall the classification of three dimensional complex Lie algebras. Fix a basis $\{e_1,e_2,e_3\}$.
The nilpotent algebra $\mathfrak{n}_3$ with the commutator matrix
$$ \left( \begin{array}{lll}
0& 0 & 1\\
0& 0 & 0 \\
0& 0 & 0  \end{array}\right).$$ 
The solvable algebra $\mathfrak{r}_{3,1}$ with the matrix
$$ \left( \begin{array}{lll}
0& 1 & 0\\
0& 0 & 1 \\
0& 0 & 0  \end{array}\right).$$ 
The simple Lie algebra $\mathfrak{sl}_2$ with the matrix
$$\left( \begin{array}{lll}
0& 0 & 1\\
0& 1 & 0 \\
1& 0 & 0  \end{array}\right).$$ 
and the projective family of pairwise non isomorphic algebras $d(r:s)$ with the matrix
$$ \left( \begin{array}{lll}
0& r & 1\\
0& 0 & s \\
0& 0 & 0  \end{array}\right).$$ 
In \cite{FP} the moduli space of these Lie algebras is described with the help of versal deformations. Let us recall the results for the nilpotent Lie algebra $\mathfrak{n}_3$. 

The cohomology spaces of the classified algebras are as follows.
\begin{equation*}
\begin{split}
&\mbox{Type}~~~~~~~~~~~~~~~~~~~~~~~~~~~H^1~~H^2~~H^3\\
&\mathfrak{n}_3~~~~~~~~~~~~~~~~~~~~~~~~~~~~~~~4~~~~5~~~~2\\
&d=\mathfrak{r}_{3,1}~~~~~~~~~~~~~~~~~~~~~~~~3~~~~3~~~~0\\
&d(1:1)=\mathfrak{r}_3~~~~~~~~~~~~~~~~~~1~~~~1~~~~0\\
&d(r:s)~~~~~~~~~~~~~~~~~~~~~~~~~1~~~~1~~~~0\\
&d(1:0)=\mathfrak{r}_2\oplus \mathbb{C}~~~~~~~~~~~~~2~~~~1~~~~0\\
&d(1:-1)=\mathfrak{r}_{3,-1}~~~~~~~~~~~~~1~~~~2~~~~1\\
&\mathfrak{sl}_2~~~~~~~~~~~~~~~~~~~~~~~~~~~~~~~~0~~~~0~~~~0
\end{split}
\end{equation*}

We get $H^2(\mathfrak{n}_3;\mathfrak{n}_3)$ is five dimensional.
Let us give explicit representative cocycles which form the basis of $H^2(\mathfrak{n}_3;\mathfrak{n}_3)$. We give the non zero values.
\begin{equation*}
 \begin{split}
&(1)~f_1: f_1(e_2,e_3)=e_3;\\
&(2)~f_2: f_2(e_1,e_2)=e_2,~ f_2(e_1,e_3)=-e_3;\\
&(3)~f_3: f_3(e_1,e_2)=e_3;\\
&(4)~f_4: f_4(e_1,e_3)=e_1;\\
&(5)~f_5: f_5(e_1,e_3)=e_2.
\end{split}
\end{equation*}

It is easy to check that all the Massey brackets are zero. So the universal infinitesimal deformation is versal  and is given by the matrix
$$\left( \begin{array}{lll}
0& t_4 & 1\\
t_2& t_5 & 0 \\
t_3& -t_2 & t_1  \end{array}\right).$$ 

Let us check that how our infinitesimal deformation (which are real deformations) fit in the moduli space of three dimensional Lie algebras.

The first deformation with cocycle $f_1$ has the matrix
$$\left( \begin{array}{lll}
0& 0 & 1\\
0& 0 & 0 \\
0& 0 & t  \end{array}\right).$$ 
which is equivalent to $\mathfrak{r}_2\oplus \mathbb{C}$.

The second deformation with cocycle $f_2$ has the matrix
$$\left( \begin{array}{lll}
0& 0 & 1\\
t& 0 & 0 \\
0& -t & 0  \end{array}\right).$$ 
which is equivalent to $\mathfrak{sl}_2$.

The third deformation with cocycle $f_3$ has the matrix

$$\left( \begin{array}{lll}
0& 0 & 1\\
0& 0 & 0 \\
t& 0 & 0  \end{array}\right).$$ 
which is equivalent to $\mathfrak{r}_{3,-1}$.

The fourth deformation with cocycle $f_4$ has the matrix
$$\left( \begin{array}{lll}
0& t & 1\\
0& 0 & 0 \\
0& 0 & 0  \end{array}\right).$$ 
which is equivalent to $\mathfrak{r}_2 \oplus \mathbb{C}$.

The fifth deformation with cocycle $f_5$ has the matrix
$$\left( \begin{array}{lll}
0& 0 & 1\\
0& t & 0 \\
0& 0 & 0  \end{array}\right).$$ 
which is equivalent to $\mathfrak{r}_{3,-1}$.

The first deformation is equivalent to $\mathfrak{sl}_2$, the second and fourth deformations give the Lie algebra $\mathfrak{r}_{3,-1}$ that means that the Lie algebra $\mathfrak{n}_3$ deforms to the family in two different ways. The third and fifth deformations are equivalent to $\mathfrak{r}_2 \oplus \mathbb{C}$. which means that $\mathfrak{n}_3$ deforms to $\mathfrak{r}_2 \oplus \mathbb{C}$ in two different ways. 

\section{Leibniz deformations of $\mathfrak{n}_3$}
The classification of three dimensional nilpotent Leibniz algebras is known. Let us recall the definition.

We take $L^1=L, L^{k+1}=[L^k, L]~\mbox{for}~k\in \mathbb{N}$.
\begin{defn}
A Leibniz algebra $L$ is called nilpotent if there exists an integer $n \in \mathbb{N}$ such that 
$$L^1 \supset L^2 \supset \ldots \supset L^n={0} .$$
The smallest integer $n$ for which $L^n=0$ is called the nilindex of $L$.
\end{defn}
The classification of complex nilpotent Leibniz algebras up to isomorphism for dimension $2$ and $3$ is in \cite{L1} and \cite{A2}.
In dimension three there are five non isomorphic algebras and one infinite family of pairwise not isomorphic algebras. The list of this classification is given below.
\begin{equation*}
 \begin{split}
\lambda_1  :  &~~\mbox{abelian}\\
\lambda_2  :  &~~[e_1,e_1]=e_2 \\
\lambda_3  :  &~~[e_2,e_3]=e_1, [e_3,e_2]=-e_1\\
\lambda_4  :  &~~[e_2,e_2]=e_1, [e_3,e_3]=\alpha e_1, [e_2,e_3]=e_1; \alpha \in \mathbb{C}\\
\lambda_5  :  &~~[e_2,e_2]=e_1, [e_3,e_2]=e_1, [e_2,e_3]=e_1\\
\lambda_6  :  &~~[e_3,e_3]=e_1, [e_1,e_3]=e_2
 \end{split}
\end{equation*}
Let us mention that in this list only $\lambda_3$ is a Lie algebra. This is the one which is denoted by $\mathfrak{n}_3$ and this is the case we are computing.

 Now we consider the Leibniz algebra $L=\lambda_3=\mathfrak{n}_3$ and compute its second Leibniz cohomology space. 

Our computation consists of the following steps:\\
 (i) To determine a basis of the space of cocycles $ZL^2(L;L)$,\\
 (ii) to find out a basis of the coboundary space $BL^2(L;L)$,\\
 (iii) to determine the quotient space  $HL^2(L;L)$.\\
(i) Let $\psi$ $\in$ $ZL^2(L;L)$. Then $\psi :L{\otimes 2}\longrightarrow L$ is a linear map and  $\delta \psi =0$, where  \begin{equation*}
 \begin{split}
 \delta \psi(e_i, e_j, e_k)
 &=[e_i,\psi(e_j, e_k)]+[\psi (e_i, e_k), e_j]-[\psi(e_i, e_j), e_k] -\psi([e_i, e_j], e_k) \\
&~ +\psi(e_i,[e_j,e_k])+\psi([e_i, e_k], e_j) ~\mbox{for}~0\leq i,j,k \leq 3.
\end{split}   
\end{equation*}
Suppose  $\psi(e_i,e_j)=\sum_{k=1} ^{3} a_{i,j}^{k} e_k$ where $a_{i,j}^{k} \in \mathbb C$
  ; for $1\leq i,j,k\leq3$.
Since $\delta \psi =0$ equating the coefficients  of $e_1, e_2
~\mbox{and}~ e_3 $ in $\delta \psi(e_i, e_j, e_k)$ we get the following relations:
\begin{equation*}
\begin{split}
&(i)~ a_{1,1}^1 =a_{1,1}^2=a_{1,1}^3=0 ;\\
&(ii)~a_{1,2}^1=-a_{2,1}^1;~a_{1,2}^2=-a_{2,1}^2;~a_{1,2}^3=-a_{2,1}^3; \\
&(iii)~a_{1,3}^1=-a_{3,1}^1; ~a_{1,3}^2=-a_{3,1}^2;~a_{1,3}^3=-a_{3,1}^3;\\
&(iv)~ a_{2,2}^2= a_{2,2}^3=0 ;\\
&(v)~a_{2,3}^2=-a_{3,2}^2; ~a_{2,3}^3=-a_{3,2}^3;\\
&(vi)~ a_{3,3}^2=a_{3,3}^3 =0 ;\\
&(vii)~a_{1,2}^2=-a_{1,3}^3
\end{split}
\end{equation*}
Observe that there
is no relation among $a_{2,2}^1$,$a_{2,3}^1$, $a_{3,2}^1$ and $a_{3,3}^1$. Therefore, in terms of the ordered basis $\{e_1\otimes e_1, e_1\otimes e_2, e_1\otimes e_3, e_2\otimes e_1, e_2\otimes e_2, e_2\otimes e_3, e_3\otimes e_1, e_3\otimes e_2, e_3\otimes e_3\}$ of $L^{\otimes 2}$ and $\{e_1, e_2, e_3\}$ of $L$,the matrix corresponding to $\psi $ is of the form
$$M= \left( \begin{array}{llrllllll}
0  & x_9 &  x_4  &-x_9 & x_6 & x_7  &-x_4   &x_{10} & x_{11} \\
0  & x_2 &  x_5  &-x_2 & 0   & x_8  &-x_5   &-x_8   & 0    \\
0  & x_3 &  -x_2 &-x_3 & 0   & x_1  &~ x_2   &-x_1   & 0   
\end{array}  \right).$$
\begin{equation*}
\begin{split}
&\mbox{where}~x_1=a_{2,3}^3;~ x_2=a_{1,2}^2 ;~ x_3=a_{1,2}^3;~ x_4=a_{1,3}^1;~
x_5=a_{1,3}^2;~ x_6=a_{2,2}^1;\\
&x_7=a_{2,3}^1;~x_8=a_{2,3}^2;~x_9=a_{1,2}^1;~x_{10}=a_{3,2}^1;~\mbox{and}~x_{11}=a_{3,3}^1
\end{split}
\end{equation*} are in $\mathbb C$~.
Let $\phi_i \in ZL^2(L;L)$ for $1 \leq i\leq 11$, be the cocycle with 
$x_i=1$ and $x_j=0$ for $i\neq j$ in the above matrix of $\psi$. It is easy to check that $\{\phi_1,\cdots, \phi_{11}\}$ forms a basis of $ZL^2(L;L)$.

(ii) 
Let $ \psi_0 \in BL^2(L;L)$. We have $\psi_0=\delta g$ for some $1$-cochain $g \in CL^1(L;L)=Hom(L;L)$. Suppose the matrix associated to $\psi_0$ is same as the above matrix $M$.

Let  $g(e_i)=g_i ^1 e_1 +g_i ^2 e_2+g_i ^3 e_3$ for $i=1,2,3$. 
The matrix associated to $g$ is given by 
\begin{center} $\left(\begin{array}{lll}
g_1 ^1&g_2 ^1  & g_3 ^1\\
g_1 ^2&g_2 ^2& g_3 ^2\\
g_1 ^3&g_2 ^3  &g_3 ^3
\end{array} \right).$ \end{center}
 From the definition of coboundary we get $$\delta g(e_i,e_j)=[e_i,g(e_j)]+[g (e_i),e_j]-\psi([e_i,e_j])$$ for $1\leq i,j \leq 3$. The matrix $\delta g$ can be written as 
 
 \begin{center}$
  \left( \begin{array}{rrrrrrrll}
0  & -g_{1}^3 & g_{1}^2  & g_{1}^3 &0  &-(g_{1}^1-g_{2}^2+g_{3}^3) & -g_{1}^2  &(g_{1}^1-g_{2}^2+g_{3}^3) &0 \\           0  & 0        & 0        & 0       &0  &-g_{1}^2                   & 0         & g_{1}^2                  &0 \\
0  & 0        & 0        & 0       &0  &-g_{1}^3                   & 0         & g_{1}^3                  &0
\end{array}  \right).$ \end{center}

Since $\psi_0=\delta g$ is also a cocycle in $CL^2(L;L)$, comparing matrices $\delta g$ and $M$ we conclude that the matrix of $\psi_0$  is of the form 
\begin{center}$
\left( \begin{array}{rrrrrrrll}
0  & x_1 & x_4 &-x_1 &0 &x_7   &-x_4 &-x_7   &0 \\
0  & 0   & 0   &0    &0 &-x_4  & 0   &x_4    &0   \\
0  & 0   & 0   &0    &0 &  x_1 & 0   &-x_1   &0
\end{array}  \right).$ \end{center}
Let ${\phi_i}^\prime \in BL^2(L;L)~\mbox{for}~i=1,4,7$ be the coboundary with $x_i=1$ and $x_j=0$ for $i\neq j$ in the above matrix of
$\psi_0$. It follows that $\{\phi_7^\prime,\phi_8^\prime,\phi_9^\prime\}$ forms a basis of the coboundary space $BL^2(L;L)$.

(iii)
It is  straightforward to check that  $$\{[\phi_1],[\phi_2],[\phi_3],[\phi_4],[\phi_5],[\phi_6],[\phi_{10}],[\phi_{11}]\}$$
span $HL^2(L;L)$ where $[\phi_i]$ denotes the cohomology class represented by the cocycle $\phi_i$. 

Thus  $dim(HL^2(L;L))=8$.

The representative cocycles of the cohomology classes forming a basis of $HL^2(L;L)$ are given explicitly below. we give the non zero values.
\begin{equation*}
 \begin{split}
&(1)~\phi_1: \phi_1(e_2,e_3)=e_3,~ \phi_1(e_3,e_2)=-e_3\\
&(2)~\phi_2: \phi_2(e_1,e_2)=e_2,~ \phi_2(e_2,e_1)=-e_2, \phi_2(e_1,e_3)=-e_3, \phi_2(e_3,e_1)=e_3;\\
&(3)~\phi_3: \phi_3(e_1,e_2)=e_3,~\phi_3(e_2,e_1)=-e_3;\\
&(4)~\phi_4: \phi_4(e_1,e_3)=e_1,~ \phi_4(e_3,e_1)=-e_1;\\
&(5)~\phi_5: \phi_5(e_1,e_3)=e_2,~ \phi_5(e_3,e_1)=-e_2;\\
&(6)~\phi_6: \phi_6(e_2,e_2)=e_1; \\
&(7)~\phi_{10}: \phi_{10}(e_3,e_2)=e_1;\\
&(8)~\phi_{11}: \phi_{11}(e_3,e_3)=e_1;\\
\end{split}
\end{equation*}

Consider, $\mu_i=\mu_0 +t \phi_i$ for $i=1,2,3,4,5,6,10,11$, where $\mu_0$ denotes the original bracket in $L$. This gives $8$ non-equivalent infinitesimal deformations of $L$.

Here $\phi_1, \phi_2, \phi_3, \phi_4, \phi_5$ are skew-symmetric, so $\phi_i \in Hom(\Lambda^2 L ; L) \subset Hom(L^{\otimes 2};L)$ for $1\leq i\leq 5$.
This are exactly the cocycles presented in the previous section. 

The last three cocycles define Leibniz deformations, more precisely their infinitesimal part. The Leibniz $2$-cocycle
$\phi_6$ defines the infinitesimal deformation with matrix 
$$\left( \begin{array}{lllllllll}
0& 0 & 0&0&t& 1&0&-1&0  \\
0& 0 & 0&0&0& 0&0&~~0&0  \\
0& 0 & 0&0&0& 0&0&~~0&0   \end{array}\right).$$ 

 $\phi_{10}$ defines the infinitesimal deformation with matrix 
$$\left( \begin{array}{lllllllll}
0& 0 & 0&0&0& 1&0&t-1&0  \\
0& 0 & 0&0&0& 0&0&~~0&0  \\
0& 0 & 0&0&0& 0&0&~~0&0   \end{array}\right).$$ 

 $\phi_{11}$ defines the infinitesimal deformation with matrix 
$$\left( \begin{array}{lllllllll}
0& 0 & 0&0&0& 1&0&-1&t  \\
0& 0 & 0&0&0& 0&0&~~0&0  \\
0& 0 & 0&0&0& 0&0&~~0&0   \end{array}\right).$$ 

It is interesting to realize that all the three Leibniz deformations are nilpotent and they are real deformation. So they can be identified with the given list. Namely the $\mu_6$ is combination of $\lambda_3$ and $\lambda_2$. $\mu_{10}$ is combination of $\lambda_2,\lambda_3$ and $\lambda_6$. $\mu_{11}$ is combination of $\lambda_2$ and $\lambda_3$.  

\begin{rem}
Notice here that $f_i=\phi_i$ for $i=1,\cdots,5$. So, if we consider the infinitesimal deformations of the Leibniz algebra $L$, that automatically contain all infinitesimal Lie algebra deformations. It is interesting to note that we get few more deformations of the original bracket $\mu_0$ giving different Leibniz algebra structures, by considering the Leibniz algebra deformation.
\end{rem} 
In order to get a simpler expression for the nontrivial cocycles let us denote $\phi_1,\cdots,\phi_6$ by $\bar{\phi}_1,\cdots,\bar{\phi}_6$, $\phi_{10}$ by $\bar{\phi}_7$ and $\phi_{11}$ by $\bar{\phi}_8$.

First we describe the universal infinitesimal deformation for $L$.
Let us denote a basis of $HL^2(L;L)^\prime$ by $\{t_i\}_{1\leq i \leq 8}$. By Remark \ref{exp of inf} the universal infinitesimal deformation of $L$ can be written as 
$$[1\otimes e_i,1\otimes e_j]_{\eta_1}=1\otimes [e_i,e_j]+ \sum_{i=1}^8 t_i \otimes \bar{\phi}_i(e_i,e_j)$$  
with base $C_1 =\mathbb{C}~\oplus_{1\leq i \leq 8} \mathbb{C}~t_i $.

\section{Extension of the infinitesimal deformation}
Let us try now to extend the universal infinitesimal deformation. 
All the $8$ infinitesimal deformations considered as one parameter deformations are real deformations as $[\bar{\phi}_i,\bar{\phi}_i]=0$ for $i=1,\cdots,8$. Now let us consider the mixed brackets $[\bar{\phi}_i,\bar{\phi}_j]$ for $i \neq j$. 
For the Lie part we get that 
$$[\bar{\phi}_1,\bar{\phi}_2],[\bar{\phi}_1,\bar{\phi}_5],[\bar{\phi}_3,\bar{\phi}_4]$$
give  nontrivial $3$-cochains from which two of them are linearly independent. We get two relations on the parameter space:
\begin{equation*}
\begin{split}
&t_1t_2+t_3t_4=0\\
&t_1t_5=0.
\end{split}
\end{equation*}

 This way the versal Lie deformation $[,]_{v_1}$ of the Lie algebra $\mathfrak{n}_3$ is defined by the infinitesimal part as follows.
\begin{equation*}
\begin{split}
&[e_1,e_2]_{v_1}=t_2e_2+t_3e_3\\
&[e_1,e_3]_{v_1}=t_4e_1 +t_5e_2-t_2e_3\\
&[e_2,e_3]_{v_1}=e_1+t_1e_3\\
\end{split}
\end{equation*}
The versal basis is the factor space $$\mathbb{C}[[t_1,t_2,t_3,t_4,t_5]]/\{t_1t_5,t_1t_2+t_3t_4\}.$$

The cocycles corresponding to only Leibniz algebras result in trivial mixed brackets:
$$[\bar{\phi}_k,\bar{\phi}_l]=0~\mbox{for}~k,l=6,7,8.$$

If we take the bracket of a cocycle from the Lie set and take the bracket with a cocycle from the Leibniz set, it turns out that not all the Massey brackets are trivial. Namely,
$$[\bar{\phi}_2,\bar{\phi}_6],[\bar{\phi}_3,\bar{\phi}_6],[\bar{\phi}_5,\bar{\phi}_6],[\bar{\phi}_3,\bar{\phi}_{10}],[\bar{\phi}_5,\bar{\phi}_{10}],[\bar{\phi}_2,\bar{\phi}_{11}],[\bar{\phi}_3,\bar{\phi}_{11}],~\mbox{and}~[\bar{\phi}_5,\bar{\phi}_{11}]$$
are nontrivial three cochains. They give us second order relation on the base of the versal deformation:
\begin{equation*}
\begin{split}
t_5t_7=0,~t_3t_6=0,~t_5t_6=0,~t_5t_8=0,~t_3t_8=0,t_2t_8,~t_3t_7~\mbox{and}~t_2t_6=0.
\end{split}
\end{equation*}
together with the relations for the Lie part we get all the second order relations for the base of the Leibniz versal deformation. 

As no higher order brackets show up we get that the versal Leibniz deformation $[,]_{v}$ is defined as follows. 
\begin{equation*}
\begin{split}
&[e_1,e_1]_v=0\\
&[e_1,e_2]_v=t_2e_2+t_3e_3\\
&[e_1,e_3]_v=t_4e_1 +t_5e_2-t_2e_3\\
&[e_2,e_1]_v=-t_2e_2-t_3e_3\\
&[e_2,e_2]_v=t_6e_1\\
&[e_2,e_3]_v=e_1+t_1e_3\\
&[e_3,e_1]_v=-t_4e_1-t_5e_2+t_2e_3\\
&[e_3,e_2]_v=(t_7-1)e_1-t_1e_3\\
&[e_3,e_3]_v=t_8e_1
\end{split}
\end{equation*}
The base of the versal Leibniz deformation is the factor space $$\mathbb{C}[[t_1,\cdots,t_8]]/<t_1t_5,t_1t_2+t_3t_4,t_2t_6,t_5t_6,t_3t_6,t_5t_7,t_3t_8,t_2t_8,t_3t_7,t_5t_8>.$$

{\bf Alice Fialowski}\\
E$\ddot{o}$tv$\ddot{o}$s Lor$\acute{a}$nd University, Budapest, Hungary.\\
e-mail: fialowsk@cs.elte.hu

{\bf Ashis Mandal}\\
Indian Statistical Institute, Kolkata, India.\\
e-mail: ashis\_r@isical.ac.in

\end{document}